\documentclass[11pt]{article}
\usepackage{amsmath,amssymb,verbatim}
\input psfig.sty

\newcommand{{\Polya}}{P\'olya{} }

\newtheorem {rema}{Remark}
\newtheorem {dfn}{Definition}

\newtheorem {lemma}{Lemma}

\newtheorem {thm}{Theorem}
\def \a   {\alpha}
\def \b   {\beta}

\def \eps {\varepsilon}

\def\E{{\mathbb{E\,}}}

\def\P{{\mathbb{P}}}
\def \R {{\mathbb{R}}}
\newcommand{\Z}{\mathbb{Z}}
\def\F{{\cal{F}}}

\def\limt{\lim_{t\to\infty}}
\def\limt0{\lim_{t\to 0}}

\def\|{\,|\,}
\def \eps {\varepsilon}
\newcommand{\BBox}{\rule{6pt}{6pt}}
\newcommand\Cox{$\hfill \BBox$ \vskip 5mm}

\def\bn#1\en{\begin{align*}#1\end{align*}}
\def\bnn#1\enn{\begin{align}#1\end{align}}

\newcommand{\ri}{\rightarrow}

\newcommand{\by}{{\bf y}}
\newcommand{\bx}{{\bf x}}

\title {Stability of a growth process  generated by monomer filling
with nearest-neighbour cooperative effects}

\author{Vadim Shcherbakov\footnote{Laboratory of Large Random Systems,
Faculty of Mechanics and Mathematics, Moscow State University, 119991, Moscow, Russia;
v.shcherbakov@mech.math.msu.su}
and Stanislav Volkov\footnote{Department of Mathematics, University of Bristol, BS8~1TW, U.K.;
S.Volkov@bristol.ac.uk}
}

\begin {document}
\maketitle

\begin {abstract}
We study stability of a growth process generated by sequential
adsorption  of particles on a one-dimensional lattice torus, that
is, the process formed by the numbers of adsorbed  particles at
lattice sites, called {\em heights}. Here the stability of
process, loosely speaking, means that its components grow at
approximately the same rate. To assess stability quantitatively,
we investigate the stochastic process formed by differences of
heights.

The model can be regarded as a variant of a \Polya urn scheme with
local geometric interaction.
\end {abstract}


\noindent {\bf Keywords:} Cooperative sequential adsorption, deposition, growth,
 urn models, reinforced random walks,  Lyapunov function.

\noindent {\bf Subject classification:} primary 60G17, 62M30;
secondary 60J20.

\section{Introduction}

\subsection{The model and results}
Let $\{1,2,\ldots,N+1\}$ be a lattice segment with periodic
boundary conditions, i.e. a one dimensional lattice torus with
$N+1$ points. Assume that $N\geq 2$. The growth process is  a
discrete-time Markov chain
$\xi(t)=(\xi_1(t),\ldots,\xi_{N+1}(t)),\, t\in
\Z_{+}=\{0,1,2,\dots\}$ with values in $\Z_{+}^{N+1}$, specified
by the following transition probabilities:
\begin{align*}
 \P\{\xi_i(t+1)=\xi_{i}(t)+1,\
\xi_j(t+1)=\xi_j(t)\ \forall j\ne i\|\xi(t)\}=
\frac{\beta^{u_i(t)}}{\sum_{j=1}^{N+1}\beta^{u_j(t)}},
\end{align*}
 \bn
 u_i(t)=\sum\limits_{j\in U_i} \xi_j(t),\ \  i=1,2,\dots,N+1,
 \en
where $\b> 0$ and  $U_i$ is  a certain neighbourhood of site $i$.

\begin{dfn}
The quantity $u_i(t)$ is called a {\em potential} of site $i$ at
time $t$.
\end{dfn}
We consider the  following three possibilities for neighbourhood $U_i$:
\begin{itemize}
\item[{\bf (A1)}]:  $U_i=\{i\}$, no interaction;
\item[{\bf (A2)}]:  $U_i=\{i, i+1\}$, asymmetric interaction;
\item[{\bf (A3)}]:  $U_i=\{i-1, i,i+1\}$, symmetric interaction,
\end{itemize}
where  $U_{N+1}=\{N+1,1\}$ due to periodic boundary conditions in
case  {\bf (A2)}; similarly $U_{1}=\{N+1, 1, 2\},\,U_{N+1}=\{N,
N+1, 1\}$   in  case {\bf (A3)}. It should be noted that periodic
boundary conditions are imposed for technical reasons only to
avoid  boundary effects.

The growth process above describes random sequential allocation of
particles at the lattice sites, where $\xi_k(t)$ is the number of
particles at site $k$ at time $t$. It is motivated by cooperative
sequential adsorption (CSA) model widely used in physics and
chemistry for representation of adsorption processes. CSA  is
probabilistic in nature and  captures the following important
feature of adsorption processes. A molecule diffusing around a
certain material surface (say, a bounded region either of
continuous space  or lattice)  might  get adsorbed by the surface.
The adsorption probability depends on a spatial configuration
formed by locations of previously adsorbed particles; for example,
the subsequent particles are more likely to get adsorbed around
locations of previously adsorbed particles. In the opposite
scenario the adsorbed particles can inhibit adsorption in their
neighbourhoods. For additional details and examples we refer the
reader to~\cite{Evans,Privman} and  references therein.

The  model under consideration relates to a version of CSA where
the unnormalized adsorption  probability at a location  depends on
the number of particles previously adsorbed in its neighbourhood.
In the mostly studied in physics adsorption model, namely, random
sequential adsorption (RSA), the adsorption probability equals $0$
at any location with one or more neighbours and equals  $1$
otherwise, i.e.\ essentially no neighbours are allowed. Asymptotic
and statistical studies    of  CSA generalizing RSA (by allowing
any number  of neighbours) were undertaken in~\cite{PenSch,Sch};
see also~\cite{Sch1}, where one proposes a model of point process
motivated by this CSA.

Our model can also be regarded as a one-dimensional lattice
variant of CSA described in~\cite{PenYuk}, where the adsorbing
probability takes the form of a product of probabilities
associated  with each of the adsorbed  nearby particles. Indeed,
in our model all neighbours  contribute the same factor $\beta>0$
to the product, so that only the number of neighbours is relevant.
If $\beta<1$, then one expects that adsorption slows down in a
saturated region.

Another  source of motivation has been  provided by monomer
filling with nearest-neighbour cooperative effects, see p.1289
of~\cite{Evans}. This is a continuous-time model on the lattice
with a hard-core type constraint: only a single particle can be
adsorbed at a site.  A site's neighbourhood  is understood as
usual, i.e.\ as in our case {\bf (A3)}, and the intensity of
adsorption at a site depends on the number of existing neighbours.
As a result, there are three non-zero intensities: $\lambda_0,
\lambda_1$ and $\lambda_2$, determining the model dynamics in the
one-dimensional case. The main difference of our model  from
monomer filling is that  we allow {\it any} number of particles to
be deposited at a site.

The infinite capacity  assumption puts our model also into the
usual ``balls and bins'' framework of urn models,
see~\cite{PemReview}, with an essential difference resulting form
additional interaction between bins. It is easy to see that in
case {\bf (A1)} the model is  a particular variant of the
well-known \Polya urn model (see a detailed discussion in Section
\ref{stab}).

Finally, we note that our model  in case {\bf (A1)} is closely
related to models of neuron growth in biology, in particular to
the one considered in~\cite{KhKh} representing the early stage of
neuron growth. In their model probability of adsorption is
proportional to the $\alpha$-th power of the number of particles
at a node ($u_i^{\alpha}(t)$ was used there instead of
$\beta^{u_i(t)}$ in our paper). The same model with polynomial
weights and its generalizations have been extensively used in
modelling  a so-called ``positive feedback'' in economics (e.g.,
see~\cite{AW-MOS,AW-O} and references therein). A special limit
variant of our  model  in case {\bf (A3)} (arising as $\beta\ri
0$, see Section~\ref{stab})  relates to the models of biological
neural networks studied in \cite{KMR,MalTur}.

\subsection{Stability}
\label{stab}

Loosely speaking,   {\it stability} of the growth process  means
that the ``profile" $\xi_i(t),\, i=1,\ldots,N+1$, is
``approximately flat", i.e. there are no extraordinary peaks. To
describe this property in a formal  way we introduce a process of
differences $\zeta(t)=(\zeta_{1}(t),\ldots,\zeta_N(t))\in \Z^N,\,
t\in \Z_{+},$ where
 $$
 \zeta_{i}(t)=\xi_i(t)-\xi_{N+1}(t),\,\, i=1,\ldots,N
 $$
and also for convenience set $\zeta_{N+1}(t)\equiv 0$. It is easy
to see that $(\zeta(t),\, t\in \Z_{+})$ is also a Markov chain
with the following transition probabilities
 \bn
 \P\{\zeta_i(t+1)=\zeta_{i}(t)+\delta_{i,k},\ i=1,2,\dots,N
 \|\zeta(t)\}=\frac{\beta^{ \sum_{j\in U_k}\zeta_j(t) }}{Z(\zeta(t))}
 \en
and
 \bn
\P\{\zeta_i(t+1)=\zeta_{i}(t)-1,\ i=1,2,\dots,N
 \|\zeta(t)\}=\frac{ \beta^{\sum_{j\in U_{N+1}}\zeta_j(t)} }{Z(\zeta(t))},
 \en
where
 \bn
 Z(\zeta(t))=\sum\limits_{k=1}^{N+1}\beta^{\sum_{j\in U_k}\zeta_j(t)}.
 \en
and
 \bn
\delta_{i,k}&=\left\{\begin{array}{rcl}
 1, &\text{if}& i=k,\\
 0, &&\text{otherwise.}
 \end{array}\right.
 \en
\begin{dfn}
We say that the growth process  is {\em stable} if the process of
differences is an ergodic (positive recurrent) Markov chain.
Otherwise the growth process is called {\em unstable}.
\end{dfn}
The following arguments apply  when $U_i=\{i\}$, when our model
becomes a particular case of generalized \Polya urn model
(see~\cite{KhKh,PemReview}). Let $\xi_i(t),\, i=1,\dots,N+1,$
represent the numbers of the balls of $N+1$ different types at
time $t$; the probability to pick a ball of a certain type $i$ is
proportional to $w(\xi_i(t))$, and here $w(x)=\beta^x$. By Rubin's
construction arguments (see \cite{BDavis}, Section~5) if
 \bn
A_i=\left\{\lim_{t\to\infty}\xi_i(t)=\infty,\  \sup_t
\xi_j(t)<\infty \ \text{ for all } j\neq i\right\},\quad
i=1,\ldots,N+1,
 \en
and $\sum_{k=1}^{\infty} w(k)^{-1}<\infty$, then
 \bn
\P\left(\bigcup_{i=1}^{N+1} A_i\right)=1.
 \en
Applying the result to  our model with $\beta>1$, since
$\zeta_i(t)=\xi_i(t)-\xi_{N+1}(t)$,  we obtain that either
$\zeta_i(t)\to \infty$ for some $i\in\{1,\dots,N\}$ (on event
$A_i$), or all $\zeta_i(t)\to-\infty$ (on event $A_{N+1}$). In
both cases transience of the process of differences follows. Thus,
the general result for the urn model implies that the growth
process is unstable.

In the opposite situation, i.e. when  $\sum_{k=1}^{\infty}
w(k)^{-1}=\infty$,  Rubin's results imply that with probability
$1$ {\it all} components of the growth process grow to infinity
and this is the case in our model with $0<\beta\leq 1$. We prove a
stronger result in the case  $0< \beta<1$, namely, we show  that
the distribution of the process of differences {\it stabilizes},
i.e. converges to a stationary distribution.

If $\beta=1$, then, regardless of the type of neighbourhood, with
probability $1$ {\it all} components of the process grow to
infinity, but the growth is unstable. Indeed, in this case the
process of differences is a zero drift spatially homogenous random
walk with bounded jumps which asymptotic behaviour is well-known
(Theorem~8.1, ch.~2, in \cite{Spitz}). If $N\geq 3$, then the
random walk is transient. If $N=2$, then the random walk is
recurrent, but it is null recurrence (as follows from the
subsequent arguments), therefore the growth is unstable for $N=2$
as well. If, just for this instance, we also allow $N=1$ (since
the process is well-defined in this case as well), then the
process of differences is just a one-dimensional simple symmetric
random walk, which is non-ergodic. If $N=2$, then null recurrence
of the random walk follows from null recurrence of its
coordinates, since each of them is just a one-dimensional
symmetric random walk. Thus in the case $\beta=1$ instability is
implied by the properties of zero-drift random walks, therefore
this case is completely eliminated from our further
considerations.

Stability of the growth  process is rather intuitive in the
no-interaction case, i.e.  $U_i=\{i\}$, if $0<\beta<1$. Indeed, in
this case growth slows down  at the sites with the maximal
potential and accelerates at the sites with the minimal potential
resulting in the stability effect, in contrast to the case
$\beta>1$ where growth accelerates  at the sites with the maximal
potential and no stability is observed. The  picture is not so
straightforward for the models with interaction: for instance, it
turns out that the growth process is unstable for the model with
symmetric interaction for any value of  $\beta$.

To understand possible sources of instability in the models with
interaction it is  helpful to consider two other growth processes,
which can be viewed as ``extreme" versions of our growth process
resulted from letting $\b\ri 0$ and $\b\ri \infty$ respectively.

An easy calculation shows that as $\b\ri 0$ the probability to get
adsorbed {\it not} at one of the minima converges to $0$.
Therefore, a natural interpretation of the formal case ``$\b=0$''
is that at time $t+1$ a particle  is allocated equally likely at
any site $i$ such that $u_i(t)=\min_{k=1,\ldots,N+1}u_k(t)$.
Similarly, ``$\b=\infty$'' can be interpreted as the situation
when at time $t+1$ a particle is allocated equally likely at any
site $i$ such that $u_i(t)=\max_{k=1,\ldots,N+1}u_k(t)$. Consider,
for instance, the case $\b=0$, $N+1=4$, and $U_i=\{i-1,i,i+1\}$.
It is easy to see that, depending on the initial configuration,
with probability $1$ one of the two following events occur: the
growth is at even nodes only, the growth is at odd nodes only.
Such limit configurations can be called an {\it attractor} of the
process by analogy with similar phenomena observed in
probabilistic models of biological neural networks
(see~\cite{KMR,MalTur} for details). If $N$ is arbitrary, then it
is also possible  to describe   {\it all} limit configurations of
the growth  processes in  both ``$\beta=0$'' and
``$\beta=\infty$'' cases. The case ``$\beta=\infty$'' is trivial
regardless of the value of $N$ and the type of interaction. On the
other hand, the limit behaviour of the growth process
corresponding to ``$\beta=0$'' is non-trivial in case of an
arbitrary $N$  for both asymmetric and symmetric interaction,
despite quite limited randomness of the process dynamics. A
detailed study of these extreme models is presented
in~\cite{SchVolkov}.

Though our study of stability of the growth process  relates also
to the study of morphology of random  interfaces of growing
materials generated by  {\it ballistic deposition} processes
(e.g., see~\cite{Penrose1,PenYuk1} and references therein), both
our setup and the methodology are quite different from the ones
used in the present paper, therefore we do not investigate  this
analogy in further details.

\subsection{Results}
\label{res} Here and further in the paper by {\it transience} of
the process $\zeta(t)=(\zeta_1(t),\dots,\zeta_N(t))$ we understand
$\lim_{t\to\infty}|\zeta(t)|\to\infty$, where $|\cdot|$ is the
usual Euclidean norm.

\begin{thm}
\label{nointeraction} Suppose $U_i=\{i\},\, i=1,\ldots,N+1$.
\begin{itemize}
\item[(1)] If $0<\beta<1$, then Markov chain $(\zeta(t),\,t\in
\Z_{+})$ is ergodic.
\item[(2)]  If $\beta>1$, then Markov chain $(\zeta(t),\, t\in
\Z_{+})$ is transient.
\end{itemize}
\end{thm}
The assertion of the second  part of Theorem \ref{nointeraction}
is  a corollary of the well known results for \Polya urn scheme,
see~\cite{AW-MOS,AW-O}, and also the discussion in
Section~\ref{stab}.

Before we formulate the next statement, we need the following
definition.
\begin{dfn}
Consider a planar process $\zeta(t)$  with polar coordinates
$({\bf r}(t),\varphi(t))$, $t\ge 0$. We say that $\zeta$ is {\rm
essentially a clockwise spiral}, if
\begin{itemize}
\item[(i)] $|{\bf  r}(t)|\to\infty$ as $t\to\infty$, and
\item[(ii)] for some large enough time $\tau_0$  such that
$\varphi(\tau_0)=2\pi k_0$ where $k_0\in\Z$ we have we have
$$
\varphi(\tau_n)=2\pi k_0+\frac{\pi n}{2}, \ \ \text{ for all }
n=0,1,2,\dots,
$$
where all
$$
\tau_n=\inf\{t>\tau_{n-1}: \varphi(t)=\frac{\pi k}{2}\text{ for
some }k\in\Z\},\ \ n=1,2,\dots
$$
are finite.
\end{itemize}
\end{dfn}

\begin{thm}
\label{asym}
Suppose $U_i=\{i,i+1\},\, i=1,\ldots,N+1$.
\begin{itemize}
\item[(1)] If $N=2$ and $0<\b<1$, then Markov chain $(\zeta(t),\, t\in
\Z_{+})$ is ergodic. Consequently, $\xi_1(t)=\xi_2(t)=\xi_3(t)$
for infinitely many $t$'s almost surely.

\item[(2)] If $\b>1$, then Markov chain $(\zeta(t),\, t\in \Z_{+})$ is transient.
Moreover, if also  $N=2$, then  the trajectory of $(\zeta_1(t),
\zeta_2(t))$ is essentially a clockwise spiral, and
$\xi_1(t)=\xi_2(t)=\xi_3(t)$ only for finitely many $t$'s a.s.
\end{itemize}
\end{thm}
It should  be  noted that when $N=2$ there is an interesting
comparison between Theorem~\ref{asym} on the one hand, and the
Friedman urn on the other hand. There will be infinitely many
``ties'' ($\xi_k(1)=\xi_k(2)=\xi_k(3)$ for infinitely many $k$'s)
if $\b<1$; in a Friedman urn with $\rho<1/2$ there will be
infinitely many ties, while the opposite occurs when $\rho>1/2$:
see \cite{Freedman} and Section 6 in \cite{MenVolk}.

\begin{thm}\label{sym} Suppose  $U_i=\{i-1,i,i+1\},\,
i=1,\ldots,N+1$. Then Markov chain $(\zeta(t),\, t\in \Z_{+})$ is
transient for any $N\geq 3$ for any $\beta\in
(0,1)\cup(1,\infty)$. Moreover, if $\beta>1$, then  with
probability $1$ there is a $k\in\{1,\dots,N+1\}$ such that
 \bn
 \lim_{t\to\infty} \xi_i(t)&=\infty, \text{ if and only if
 }i\in\{k-1,k\}, \mbox{ and}\\
 \lim_{t\to\infty} \frac{\xi_k(t)}{\xi_{k-1}(t)}&=\beta^{c},
 \en
where $c=\lim_{t\to\infty} [\xi_{k+1}(t)-\xi_{k-2}(t)]\in \Z$.
\end{thm}

To prove the results of the present paper we combine  the
constructive methods of studying asymptotic behaviour of countable
Markov chains from~\cite{fmm} with probabilistic techniques used
in the theory of processes with reinforcement from~\cite{V2001}
(in contrast with the purely combinatorial methods used
in~\cite{SchVolkov}).

\section{Proofs}

\subsection{Proof of Theorem \ref{nointeraction}}

If $U_i=\{i\}$, then process of differences $(\zeta_t,\, t\in
\Z_{+})$ has  the transition probabilities
$$
\P\{\zeta_i(t+1)=\zeta_{i}(t)+\delta_{i,k},\ \forall  i
\|\zeta(t)\}=\frac{\beta^{\zeta_k(t)}}{1+\sum_{j=1}^{N}\beta^{\zeta_j(t)}},
\quad k=1,\ldots,N,
$$
and
$$
\P\{\zeta_i(t+1)=\zeta_{i}(t)-1,\,
i=1,\ldots,N\|\zeta(t)\}=\frac{1}{1+\sum_{j=1}^{N}\beta^{\zeta_j(t)}}
$$
Suppose $0<\beta<1$.  Consider the following function
$$
f(\bx)=|\bx|^2=\sum\limits_{i=1}^{N}x_i^2, \quad \bx=(x_1,\ldots,x_N)\in \Z^N.
$$
It is easy to see that
\begin{equation*}
 \Delta:=\E[f(\zeta(t+1))-f(\zeta(t))\| \zeta_t=\bx]=\frac{\sum_{i=1}^{N}
 \left(2x_i (\beta^{x_i} -1)+1+\beta^{x_i}\right)} {1+\sum_{i=1}^N
\beta^{x_i} },
\end{equation*}
and for any $\eps>0$
\begin{equation*}
 \Delta+\eps=\frac{\eps-\sum_{i=1}^{N}
(2x_i (1-\beta^{x_i})-1-(\eps+1)\beta^{x_i})} {1+\sum_{i=1}^N
\beta^{x_i} }.
\end{equation*}
Now let
$$
 h(x)=2x(1-\beta^x)-1-(\eps+1)\beta^x.
$$
It is clear that for $x>0$ sufficiently large, $h(x)\sim 2x$,
formally, there is an $A'>0$ such that for $x\ge A'$ we have
$h(x)\ge x$. Conversely, when $x=-a<0$,
$$
  h(x)=h(-a)=\frac 1{\beta^a} (2a-\eps-1)-1-2a.
$$
so that $h(-a)$ grows approximately exponentially and again there
is an $A''>0$ such that $h(-a)\ge a$ when $a\ge A''$. Now set
$$
 -C:=\inf_{x\in(-A',A')} (h(x)-|x|)=\inf_{x\in\R}  (h(x)-|x|)
$$
(note that $\eps+2\le C<\infty$). Consequently,
 \bn
 \sum_{i=1}^{N} (2x_i (1-\beta^{x_i})-1-(\eps+1)\beta^{x_i})&=\sum_{i=1}^{N} h(x_i)
 =\sum_{i=1}^{N} |x_i|+\sum_{i=1}^{N} \left(h(x_i)-|x_i|\right)
 \\
 &\ge \sum_{i=1}^{N} |x_i| -CN
 \en
Therefore,
\begin{eqnarray*}
 \Delta+\eps\le\frac{CN+\eps-\sum_{i=1}^{N}
|x_i|} {1+\sum_{i=1}^N \beta^{x_i} }\le 0
\end{eqnarray*}
except for a possibly finite number of $(x_1,\dots,x_N)$ lying in
the ``bad'' set
$$
M:= \left\{\sum_{i=1}^{N}|x_i|\le CN+\eps\right\}.
$$
Hence, the conditions of the Foster criterion
(Theorem~\ref{th_fa1})  are satisfied and Markov chain $(\zeta(t),\, t\in \Z_{+})$ is ergodic.
Thus the first part of Theorem \ref{nointeraction} is now proved.
  \Cox

\subsection{Proof of Theorem \ref{asym}}
\label{sec_proof_asym}

\subsubsection{Proof of part (1) of Theorem \ref{asym}}

\begin{figure}[htb]
\vspace{1cm}
\centerline{\hbox{\psfig{figure=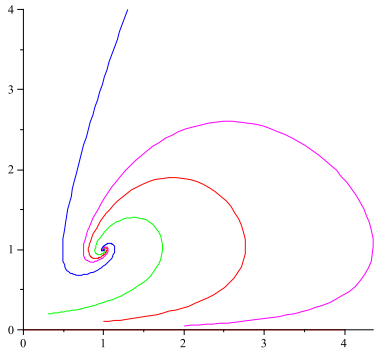,width=5cm,height=5cm}}}
\caption{Trajectories of the dynamical process following the
Markov chain.} \label{Fig2}
\end{figure}

Let $N=2$ and $0<\beta<1$. Set $\zeta_1(t)=\xi_1(t)-\xi_3(t)$ and
$\zeta_2(t)=\xi_2(t)-\xi_3(t)$. The new process
$\zeta(t)=(\zeta_1(t), \zeta_2(t))$ is a time-homogeneous Markov
chain on $\Z^2$ with the following transitions
 \begin{align*}
\P\{\zeta(t)=(x_1+1, x_2)|\zeta(t)=(x_1,x_2)\} & = \frac{\beta^{x_1+x_2}}{\beta^{x_1+x_2}+\beta^{x_2}+\beta^{x_1}},\\
\P\{\zeta(t)=(x_1, x_2+1)|\zeta(t)=(x_1,x_2)\}& =
 \frac{\beta^{x_2}}  {\beta^{x_1+x_2}+\beta^{x_2}+\beta^{x_1}}\\
\P\{\zeta(t)=(x_1-1, x_2-1)|\zeta(t)=(x_1,x_2)\} & =
\frac{\beta^{x_1}}  {\beta^{x_1+x_2}+\beta^{x_2}+\beta^{x_1}}.
 \end{align*}
It is natural to expect that the Markov chain approximately
follows the solutions of the differential equation
 \bn
 \frac{dx_2}{dx_1}=\frac{\b^{x_2}-\b^{x_1}}{\b^{x_1+x_2}-\b^{x_1}}
 \en
or, after making a substitution $u=\b^{-x_1}$, $v=\b^{-x_2}$,
 \begin{align}\label{eq_dyn_sys}
 \frac{dv}{du}=\frac{v(u-v)}{u(1-v)}
 \end{align}
which we cannot solve analytically, but whose solutions seem to be
spirals (see Figure~\ref{Fig2}). Hence a good candidate for a
Lyapunov function for the process, i.e.\ the function
$f(\bx)=g(\b^{-x_1},\b^{-x_2})$ such that
$\E[f(\zeta(t+1))-f(\zeta(t))\|\zeta(t)=\bx] \le -{\rm const}$
would be a function which level curves have a constant angle with
the vector field generated by (\ref{eq_dyn_sys}). Then the level
curves for $g(u,v)$ will satisfy the differential equation

\begin{figure}[htb]
\vspace{1cm}
\centerline{\hbox{\psfig{figure=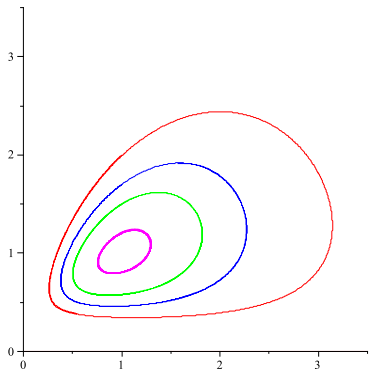,width=5cm,height=5cm}}}
\caption{``Ideal'' level curves.} \label{Fig3}
\end{figure}

\begin{figure}[htb]
\vspace{0cm}
\centerline{\hbox{\psfig{figure=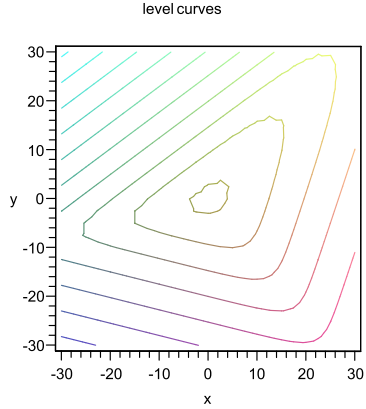,width=6cm,height=6cm}}}
\caption{Level curves for $f(x,y)$ given in (\ref{eq_ergod}) for
$\b=0.9$} \label{Fig1}
\end{figure}

 \bn
  \left(\frac{\partial g}{\partial u} , \frac{\partial g}{\partial
  v}\right)
  \cdot \left[ (v(v-u), u(1-v))
 \left(\begin{array}{cc} \cos \a & -\sin\a   \\ \sin \a & \cos\a  \end{array}\right)
 \right]=0
 \en
for some $\a$, which in turn might depend on initial conditions.
Numerical solutions for the level curves are presented in
Figure~\ref{Fig3}.  Though not being able to solve the above
equations analytically, we found  an alternative suitable function
(\ref{eq_ergod}), whose level curves are in Figure~\ref{Fig1}.
Thus the proof of part (1) of Theorem~\ref{asym} is based on the
following lemma.

\begin{lemma}
\label{lem_b1}
Suppose that $\beta<1$. Let
\begin{equation}\label{eq_ergod}
f(\bx)=\b^{1-x_1-x_2}+\b^{-3x_1+x_2}+\b^{3x_1-4x_2}+\b^{x_1+4x_2},
\end{equation}
for any  $\bx=(x_1,x_2)$.
Then for any $\delta>0$
$$
 \E[ f(\zeta(t+1))-f(\zeta(t))\|\zeta(t)=\bx] \le -(1-\b)
$$
once $|\bx|\equiv \sqrt{x_1^2+x_2^2}$ is larger than some
$C=C(\b)>0$.
\end{lemma}

{\sf Proof of Lemma~\ref{lem_b1}.} We see that
 \begin{align}\label{eq_snos}
 \E[ f(\zeta(t+1))-f(\zeta(t))\|\zeta(t)=\bx]+1-\b=-\frac{(1-\b)A(u,v)}{ \b^5 v^4 u^3 (1+u+uv)}
 \end{align}
where
 \bn
 A(u,v)&=  (\b^7+\b^6)v^3 u^3 + (\b^6+\b^5 ) v^5u +( \b^8+\b^5+\b^6+\b^7 ) v^9u^4 \\
 & + [(\b^7+\b^5+\b^6) u- (\b+\b^2+\b^3+\b^4 )] v u^6 +v^6\b^5+\b^5u^7 \\
 &   + \b^5 v^9 u^5 - \b^5 v^4 u^3 \\
 &-(1+u^2+uv+u^2v)u^2 v^4 \b^5- (1+\b+\b^2+\b^3+\b^4)u^5{v}^8  - (\b^2+\b^4+\b^3) uv^6
 \en
and  $u=\b^{x_1}$, $v=\b^{x_2}$. The term $\b^5 u^5 v^9$ clearly
dominates all the negative terms, and the term in the square
brackets is nonnegative for $u$ sufficiently large,  hence there
is a constant $c_1$ such that $A(u,v)>0$ once $\min(u,v)\ge c_1$.

Similarly,
 \bn
 u^7 A(1/u,v)&= \b^5 u^7 v^6 +O(u^6 v^6)+  (\b^8+\b^5+\b^6+\b^7)v^9 u^3 + O(v^9 u^2)
 \en
which is also nonnegative once $\min(u,v)\ge c_1$ for some $c_2>0$
and
 \bn
 v^9 A(u,1/v)&= \b^5 u^7 v^9 +O(u^7 v^8)
 \en
is  nonnegative when $\min(u,v)\ge c_3$ for some $c_3>0$.

Finally,
 \begin{align}\label{eq1u1v}
 u^7 v^9 A(1/u,1/v)&= \b^6 u^6 v^4 + \b^7  u^4 v^6 +(\b^5+\b^6+\b^7)v^8 \\
 &+ \b^5 v^4( v^5 +u^6 -[\b^{-4}+\b^{-3}+\b^{-2}+\b^{-1}] v^4 u ) \nonumber \\
 &+ \b^5 u^4 v^3 (\b v^3+u^3-v^2 u)+ O(u^6 v^3)+O(u^4 v^5) \nonumber
 \end{align}
Let $K=\b^{-4}+\b^{-3}+\b^{-2}+\b^{-1}$. For $v\ge K u$ the
expression in the second line of (\ref{eq1u1v}) is always
nonnegative; on the other hand for $u> K^{-1} v$ the term $u^6$ is
going to dominate in this line, hence for $u,v$ larger than $K^6$
the second line of (\ref{eq1u1v}) is nonnegative. The third line
of (\ref{eq1u1v}) is nonnegative when $u\ge v$. In principle, it
{\it can} become negative when $u<v$, however, since $x$ and $y$
are integers, it implies that when $u<v$, also $u\le \b v$, so
that $\b v^3+u^3-v^2 u> v^2(\b v -u)\ge 0$. Thus we have
established that $A(1/u,1/v)\ge 0$ for all
legitimate\footnote{i.e.\ of the form $u=\b^{x_1}$, $v=\b^{x_2}$,
$(x_1,x_2)\in\Z^2$.} $u,v$ such that $\min\{u,v\}\ge c_4$.

Consequently, we have shown that on the positive quadrant
$A(u,v)\ge 0$ whenever
$$
(u,v)\in [\eps^{-1},\infty)\times [\eps^{-1},\infty)\
 \cup \ (0,\eps]\times (0,\eps] \
 \cup \ (0,\eps]\times [\eps^{-1},\infty)\
 \cup \ (\eps^{-1},\infty)\times (0,\eps]
$$
for some $\eps>0$ (assume without loss of generality that
$\eps<1$).

Next, since
$$
A(u,v)=(v \b^5+\b^5+\b^7 v+\b^6 v) u^7 +O(u^6)
$$
as a Taylor series on $u$ and
$$
A(u,v)=u^4 \b^5 (\b^3+\b^2+\b+u+1) v^9 +O(v^8)
$$
as a Taylor series on $v$, we see that there is an
$\eps_1\in(0,\eps]$ such that for all $v\in [\eps,\eps^{-1}]$ and
$u\ge \eps_1^{-1}$, and for all $u\in [\eps,\eps^{-1}]$ and $v\ge
\eps_1^{-1}$ respectively, we have $A(u,v)\ge 0$.

On the other hand, for small $u$ and $v$ respectively, we have
 \bn
  A(u,v)&= \b^5 v^6+ u \times {\sf Polynom}_1(u,v,\b) =\b^5 u^7+v\times {\sf
  Polynom}_2(u,v,\b).
 \en
where ${\sf Polynom}_i(u,v,\b)$, $i=1,2$, are some polynomial
expressions involving $u$, $v$, and $\b$. Hence there is an
$\eps_2\in(0,\eps]$ such that $A(u,v)\ge 0$ whenever $0<u\le
\eps_2$ and $v\in [\eps,\eps^{-1}]$ or $0<v\le \eps_2$ and $u\in
[\eps,\eps^{-1}]$.

Combining all the results above, we conclude that the right side
of (\ref{eq_snos}) is nonnegative for all $\bx\in\Z^2$ which are
outside of the square $\left[-R,R\right]^2$ where $ R=\log(
\min\{\eps_1,\eps_2\})/\log(\b). $ Lemma~\ref{lem_b1} is proved.
\Cox

Because of Lemma~\ref{lem_b1} and the fact that $f(\bx)\to \infty$ whenever
$|\bx|\to\infty$, we can apply Foster's supermartingale criterion (Theorem~\ref{th_fa1}). \Cox

\begin{figure}[htbp]
\vspace{1cm}
\centerline{\hbox{\psfig{figure=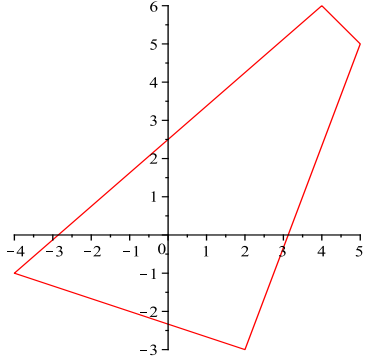,width=5cm,height=5cm}}}
\caption{Level curve for $\tilde f(x,y)$ given in
(\ref{eq_erga}).} \label{Fig4ugla}
\end{figure}

\begin{rema}We observed that the following function
 \begin{align}\label{eq_erga}
\tilde f(\bx)=\max\left\{\frac{x_1+x_2}{10},
\frac{8x_1-3x_2}{25},\frac{8x_1-7x_2}{20},\frac{-x_1-3x_2}{7}\right\}
 \end{align}
also satisfies Foster's  supermartingale criterion (its level
curves are congruent to the ones in Figure~\ref{Fig4ugla});
however, the proof of this fact requires going through a large
number of special cases and hence is omitted in favour of using
(\ref{eq_ergod}).
\end{rema}

\subsection{Proof of part (2) of Theorem~\ref{asym}}
 We start by proving
transience for any $N\ge 2$. Set
$\eta_i(t)=\xi_i(t)-\xi_{i-1}(t),\, i=2,\ldots, N+1$ and
$\eta_{1}(t)=\xi_1(t)-\xi_{N+1}(t)$, because of the periodic
boundary conditions. Then
 \bn
\xi_1(t)&=\xi_{N+1}(t)+\eta_1(t), \\
\xi_2(t)&=\xi_{N+1}(t)+\eta_1(t)+\eta_2(t), \\
 \dots &= \dots, \\
 \xi_{N+1}(t)&=\xi_{N+1}(t)+\eta_1(t)+\dots+\eta_{N+1}(t),
 \en
and obviously $\sum_{i=1}^{N+1} \eta_i(t)=0$.  We also have
$$
\P\{\xi_i(t+1)=\xi_i(t)+1|\xi(t)\}\propto
 \b^{2[\eta_1(t)+\dots+\eta_i(t)]+\eta_{i+1}(t)}
$$
for $i\le N$, and, in turn,
 \bn
 \xi_i(t+1)= \xi_i(t)+1 \Leftrightarrow
 \left\{\begin{array}{rcl}\eta_i(t+1)&=& \eta_i(t)+1,\\
 \eta_{i+1}(t+1)&=& \eta_{i+1}(t)-1.\end{array}\right\}
 \en
We are going to show that Markov chain
$$
\eta(t)=\left(\eta_1(t),\eta_2(t),\dots,\eta_{N+1}(t)\right), \quad t\in \Z_{+},
$$
with state space  $\{\bx=(x_1,\dots,x_{N+1}):\
x_1+\dots+x_{N+1}=0\}$ is transient.
It is easy to see that transience of $(\eta(t), t\in \Z_{+})$
 implies transience of $(\zeta(t), t\in \Z_{+})$.
Indeed, if $(\zeta(t), t\in \Z_{+})$ were recurrent, then
$\xi_1(t)=\dots=\xi_{N+1}(t)$ for infinitely many $t$ almost
surely. But this would imply that
$\eta_1(t)=\dots=\eta_{N+1}(t)=0$ for infinitely many $t$ almost
surely as well, thus contradicting its transience. In other words,
coordinates  of $\eta$-process are consecutive  differences of
coordinates of the growth process, therefore transience of
$\eta(t)$ yields transience of  the process of differences,
because the probability distribution of the later does not depend
on the subtracted coordinates (due to the symmetry and periodic
boundary conditions).

To prove transience of $(\eta(t), t\in \Z_{+})$ we are going  to apply Theorem~\ref{th_fa2}. To this end,
consider the function
 \bn
f(\bx)=\left\{\begin{array}{ll}
 1, &\text{ if }x_1=x_2=\dots=x_{N+1}=0\\ \\
 \frac {\displaystyle 1}{\displaystyle \max\{|x_1|,|x_2|,\dots,|x_{N+1}|\}}, &\text{ otherwise.}
\end{array}\right.
 \en
defined for any $\bx=(x_1,\ldots, x_{N+1})$. We now will show
that, provided that $|\bx|$ is large enough, $\E
[f(\eta(t+1))-f(\eta(t))\|\eta(t)=\bx]\le 0$, establishing the
result.

Indeed, consider the four following possibilities.
\begin{itemize}
\item[(a1)] $N\geq 3$ and the maximum $L$ of $x_1, x_2,\dots, x_{N+1}$ is not
unique  and  there are at least two indices $i,j$ such
that $|x_i|=|x_j|=L$ and they are at least distance $2$ apart,
i.e., $|i-j|\ge 2$ (understood with periodic boundary conditions).
\item[(a2)] $N=2$ and $|x_1|=|x_2|=|x_3|=L$.
\item[(b)]  The maximum $L$ is not unique and there are exactly two
maximums next to each other.
\item[(c)] The maximum $L$ is unique.
\end{itemize}

Observe that in cases (a1) and (a2)  no step of the chain $\eta$
can decrease the maximum, hence $\E[f(\eta(t+1)) -
f(\eta(t))\|\eta(t)=\bx]\le 0$ here.

In case (b) suppose without loss of generality that the maximum is
achieved at nodes $2$ and $3$, so $|x_2|=|x_3|=L$. Then
$\E[f(\eta(t+1)) - f(\eta(t))\|\eta(t)=\bx]$ is proportional (up
to a positive coefficient) to
 \bn
 & \b^{2x_1+x_2} \left[\frac1{\max\{|x_1+1|,|x_2-1|,L\}}-\frac 1L\right]
 \\&+
 \b^{2x_1+2x_2+x_3} \left[\frac1{\max\{|x_2+1|,|x_3-1|\}}-\frac 1L\right]
\\&+
 \b^{2x_1+2x_2+2x_3+x_4} \left[\frac1{\max\{L,|x_3+1|,|x_4-1|\}}-\frac
 1L\right]
 \en
(if $N=2$ then we set $x_4\equiv x_1$). We can see that unless
$x_3=L=-x_2$ the quantity above is always negative. However, when
$x_3=L=-x_2$ we also have that the quantity above is proportional
to
 \bn
 &\b^{-x_2-x_3} \left[\frac1{L+1}-\frac 1L\right]
+
 \left[\frac1{L-1}-\frac 1L\right]
+
 \b^{x_3+x_4} \left[\frac1{L+1}-\frac1L\right]\\
&= \frac {2-(L-1)\b^{x_3+x_4}}{L(L^2-1)}\le \frac {2-(L-1)\b}
{L(L^2-1)}
 \  \text{\ (since $|x_4|\le L-1\ \Rightarrow\ x_3+x_4\ge 1$)}
 \en
which is negative once $L\ge 3$.

Finally, in case (c) suppose that the maximum is achieved at node
$2$, $|x_2|=L$. Then $\E[f(\eta(t+1)) - f(\eta(t))\|\eta(t)=\bx]$
is proportional to
 \bn
  (*):=& \left[\frac1{\max\{|x_1+1|,|x_2-1|\}}-\frac 1L\right]
+
 \b^{x_2+x_3} \left[\frac1{\max\{|x_2+1|,|x_3-1|\}}-\frac 1L\right]
 \en
When $x_2=L$,
 \bn
 (*)&= \left[\frac1{\max\{|x_1+1|, L-1\}}-\frac 1L\right]
 +
 \b^{x_2+x_3} \left[\frac1{L+1}-\frac 1L\right]
 \\
 &\le \frac 1{L(L-1)}- \frac{ \b }{L(L+1)}  \le 0
 \en
once $L>(\b+1)/(\b-1)$.

And in the last subcase $x_2=-L$
 \bn
 (*)&= \left[\frac1{L+1}-\frac 1L\right]
 +
 \b^{x_2+x_3} \left[\frac1{\max\{L-1,|x_3-1|\}}-\frac 1L\right]
 \\
 &\le -\frac 1{L(L+1)}+ \frac{ \b^{-1} }{L(L-1)}  \le 0
 \en
again once $L>(\b+1)/(\b-1)$.

Consequently, $\E[f(\eta(t+1))\|\eta(t)=\bx]\le f(\bx)$ whenever
$|\bx|$ is sufficiently large, and by Theorem~\ref{th_fa2} Markov
chain $(\eta(t),\, t\in \Z_{+})$ is transient. So the transience
is proved for any $N\geq 2$.

\vskip 5mm

From now on assume that $N=2$. First, let us prove that the
process $\zeta(t)=(\zeta_1(t), \zeta_2(t))$ cannot remain
indefinitely in either of the following $6$ areas: $x>y>0$,
$y>x>0$, $x>0>y$, $y>0>x$,  $0>y>x$, and $0>x>y$.

Indeed, suppose that for some time $t_0$  we have $\zeta(t_0)\in
\Psi_1:=\{(x,y)\in \R^2: x>y>0\}$. When $\zeta(t)\in\Psi_1$, it is
clear that $\zeta_2(t)$ has the property
$\E[\zeta_2(t+1)\|\zeta(t)=(x,y)]\le \zeta_2(t)$. Hence
$M(t)=\zeta_2(t\wedge \tau)$, $t\ge t_0$, where
$$
\tau:=\inf\{t>t_0:\ \zeta_1(t)=\zeta_2(t)\text{ or }
\zeta_2(t)=0\}=\inf\{t>t_0: \zeta(t)\notin \Psi_1\}
$$
is a non-negative supermartingale which converges a.s. Since
$M(t+1)- M(t)$ takes only integer values, it means that for some
(random) $T$, $M(t)=const$ for all $t\ge T$. However, this is
impossible unless $\tau<\infty$, since for a fixed $y>0$ the
probability
$$
\P(M(t+1)=M(t)-1\| \zeta(t)=(x,y)\in\Psi_1)= \frac
1{1+\beta^y+\beta^{y-x}}>\frac 1{2+\beta^y}
$$
does not go to zero, thus implying by the conditional
Borel-Cantelli lemma (see e.g.~Corollary 2 on p.~518
in~\cite{Shyr}) that for some large $t>T$ we have $M(t)\ne
M(t+1)$. Therefore $\tau$ is finite.

When $\zeta(t_0)\in \Psi_2:=\{(x,y)\in \R^2: y\ge x\ge 0\}$, set
$M(t)=\zeta_2(\tau\wedge t)-\zeta_1(\tau\wedge t)$ and
$$
\tau:=\inf\{t>t_0:\ \zeta_1(t)=\zeta_2(t)\text{ or }
\zeta_1(t)=0\}.
$$
Again $M(t)$ is a non-negative supermartingale  which cannot
converge unless $\tau<\infty$, since on $\{\tau=\infty\}$ the
event $|M(t+1)-M(t)|<1$ implies $M(t+1)=M(k)$ and thus
$\zeta_1(t+1)=\zeta_1(t)-1$.

The remaining $4$ cases will analyzed briefly, since the argument
is very similar. When $x>0>y$ or $y>0>x$, the probability to go
towards the axis $x=0$ is larger than away from it, hence
eventually $\zeta(t)$ will leave either of these areas. When
$0>y>x$ there is a drift upwards, and finally when $0>x>y$ there
is a drift towards the line $x=y<0$ so that a non-negative
supermartingale $\zeta_1(t)-\zeta_1(t)$ must converge, and at the
same time $\zeta(t)$ cannot remain indefinitely on the line
$(-a,-a-\Delta)$, where $\Delta>0$ is a constant and
$a=i,i+1,i+2,\dots$, since along this line there is a
non-diminishing probability to go up, of order $\b^{-\Delta}$.

Finally, we need to show that the process $\zeta(t)$ eventually
``rotates'' in one direction on its way to infinity. This
immediately follows from the following Lemma~\ref{lem_rot}, taking
into account the established transience, the Borel-Cantelli lemma
and the summability in $a$ of the terms on the RHS
of~(\ref{eq_sumable}). \Cox

\begin{lemma}\label{lem_rot}
Assume  $\b>1$. Suppose $a,k\ge 1$ and let $E_k$ be the event
$\{$the trajectory of $\zeta(t)$ for $t>k$ crosses $\{(0,y),y<
0\}$, $\{(x,0),x<0\}$, $\{(0,y),y> 0\}$, $\{(x,0),x\ge 0\}$
exactly in this order and ends up at point $(a',0)$ with $a'\ge
a+1\}$. Then, for some $\nu=\nu(\beta)>0$ and all large $a>0$ we
have
 \begin{align}\label{eq_sumable}
 \P\left(E_k^c\| \zeta(k)=(a,0)\right)\  \le \exp\left(-\nu a\right).
 \end{align}
Also, similar statements hold for starting points $(0,-a)$,
$(-a,0)$ and $(0,a)$.
\end{lemma}
{\sf Proof.} (1) Let us show that when the process leaves $(a,0)$,
it ends up at $(0,-b)$, $b\ge a$, with probability very close to
$1$. This will be done by demonstrating that initially the process
moves only right and left-down (south-west, SW for short), and
then after reaching level $y=-a/2$ (for simplicity assume that $a$
is even), it moves only in the SW direction, all with high
probability.

Indeed, suppose that $\zeta(\tau_0)=(a,0)$ for some $\tau_0$. For
$i=1,2,\dots$ define
 \bn
 \tau_{i}&=\inf\{t:\ \zeta_2(t)=-i\},
 \\
 A_{i}&=\left\{\exists k\in\{1,\dots,a\}:\
 \zeta_2(\tau_{i-1}+k)=-i,\text{ but }
 \zeta_2(\tau_{i-1}+m)=-i+1\ \forall m\in\{0,\dots,k-1\}\right\},
 \en
where the latter is the event that after reaching level
$y=-(i-1)$, the process $\zeta(t)$ makes less than $a$ consecutive
steps to the right, after which it moves in the SW direction.

On the event $A_{i}$ the stopping time $\tau_{i}$ is finite, also
on $\cap_{j=1}^{i} A_j$  we have
$$
 \zeta(\tau_{i})\in G:=\{(x,y)\in \Z^2: y\le 0\text{ and } y\le x-a\le
 -ay\}.
$$
Then for $1\le i\le a/2-1$
 \bnn\label{eq_A|A}
 \P\left(A_{i+1}\| \bigcap_{j=1}^{i} A_j\right)\ge
  \left(1-\b^{-ia}\right)
  \left(1-\b^{-a}\right)\ge (1-\b^{-a})^2
 \enn
since the probability to jump to the right from $(x,-i)$ is
$$
 \frac{\b^{x-i}}{\b^{-i}+\b^x+\b^{x-i}}<\b^{-i}.
$$
and the conditional probability to jump up from $(x,-i)$, given
that the next jump is {\it not} rightwards, is
$$
 \frac{\b^{-i}}{\b^{-i}+\b^x}<\b^{-x-i}\le \b^{-a}
$$
(recall that  $\zeta(\tau_{i})\in G$). Taking into account that
$$
\P(A_1)\ge \left(1-\frac 1{2^a}\right)\left(1-\b^{-a}\right)
$$
from (\ref{eq_A|A}) we conclude
 \bnn\label{eq_L21}
 \P\left( \bigcap_{j=1}^{a/2} A_j\right)\ge 1- a e^{-
 \tilde\nu a}
 \enn
where
$$
\tilde \nu=\min\left\{\log \b, \log 2\right\}>0.
$$

Once the process $\zeta(t)$ has reached level $y=-a/2$ such that
all $A_i$, $i=1,2,\dots,a/2$, occurred, we know that
$\zeta(\tau_{a/2})=(x^*,-a/2)$ and $a/2\le  x^* \le a+a^2/2$.
However, from $(x,-i)$ the probability {\it not}\/ to go SW is
$$
 \frac{\b^{-i}+\b^{x-i}}{\b^{-i}+\b^x+\b^{x-i}}\le \b^{-i}
$$
(as $x\ge i$) therefore with probability at least
 \bnn\label{eq_L22}
 \left(1-\b^{-a/2}\right)^{x^*}&\ge  \left(1-\b^{-a/2}\right)^{a+a^2/2}\ge 1-(a+a^2/2)\b^{-a/2}
 \ge 1-a^2 e^{-\tilde\nu a/2}
 \enn
the process will make only SW steps until it reaches the vertical
axes at some point $(0,-b)$ with $b\ge a$.

The argument for the remaining 5 cases is very similar, hence we
just sketch the proof.
 \vskip 1mm

 (2) Suppose the process starts at $(0,-a)$, $a\gg 1$.
We will show that with probability close to $1$ it ends up on the
line $x=y$ at the point $(-b,-b)$ with $b\ge a+1$.

Consider the trajectory of the process along the lines
$(-i,-i-j)$, where $j$ is called a {\it level}, with the
probabilities to jump up, right, and SW being proportional to
$\b^{-j},\b^{-i-j},1$, respectively. Also, initially $i=0$, $j=a$.
With probability
$$
(1-2\b^{-a})^a\ge 1-2a\b^{-a}
$$
the process starts with consecutive $a$ SW steps. After this, a
conditional probability that the jump was ``right'' but not
``up'', given that one of the two has indeed occurred, is at most
$\b^{-a}$. Hence, with probability
$$
(1-\b^{-a})^a\ge 1-a\b^{-a}
$$
the process will consecutively pass through the levels $j-1$,
$j-2$, $\dots$, $0$ until it reaches the line $x=y<0$. Is is also
clear that the process cannot remain indefinitely on the same
level $j$, as there is a constant probability to go up, of order
$\b^{-j}$. Consequently, the process will arrive to $(-b,-b)$
where $b\ge 2a$ (in fact, we expect $b$ to be of order
$1+\b+\dots+\b^a$).
 \vskip 1mm

(3) Next, suppose the process starts at $(-a,-a)$, $a\gg 1$. We
will show that with probability close to $1$ it ends up on the
line $y=0$ at the point $(-b,0)$ with $b\ge a$.
 \vskip 1mm

Again, notice that after approximately geometrically $(1/2)$
distributed number of SW steps along the line $x=y$, the process
jumps up and the relative probability to jump ``up'' vs.\ ``SW''
is around $\b/(\b+1)>1/2$. As before, we can show that the process
will reach the horizontal axis in a number of steps of order $a$
never ever making a ``right'' step with probability of exactly the
same order as in (1).
 \vskip 1mm

(4) Suppose the process starts at $(-a,0)$, $a\gg 1$. By very
similar arguments we can show that the process will reach the
point $(0,b)$ with $b\ge a$ never making a SW move on its way,
thus going through level lines $x=-j$, $k=a,a-1,\dots,0$, with
probability close to $1$. It should take around
$1+\b+\dots+\b^{a}$ steps.
 \vskip 1mm

(5)  Suppose the process starts at $(0,a)$, $a\gg 1$. Again, by
similar arguments one can show that the process will reach the
point $(b,b)$ with $b\ge a$ never making a SW move on its way,
with probability close to $1$.
 \vskip 1mm

 (6)  Finally, suppose
the process starts at $(a,a)$, $a\gg 1$. One can easily show that
the process will reach the point $(b,0)$ with $b\ge a$ never
making an ``up'' move on its way, with probability close to $1$.
\vskip 1mm

To finish the proof, fix some $\nu\in(0,\tilde\nu/2)$. Combining
(\ref{eq_L21}) and~(\ref{eq_L22}) with similar results established
in (2) through (6), we conclude that the probability that
$\zeta(t)$ started at $(a,0)$ will sequentially visit the areas
$\{x>0,y<0\}$, $\{x<0,y<0\}$, $\{x<0,y>0\}$, $\{x>0,y>0\}$, and
end up at a point $(a',0)$ with $a'\ge a+1$ is at least
$$
1- e^{-\nu a}
$$
for all $a$ larger than some constant $A=A(\nu,\b)$.\Cox

\subsection{Proof of Theorem \ref{sym}}

\subsubsection{Proof of Theorem \ref{sym} for $\beta>1$}
Recall that in the  symmetric case
$$
\P\{\xi_i(t+1)=\xi_{i}(t)+ \delta_{i,k},\ \forall  i \|
\xi(t)=\xi\}\varpropto \beta^{\xi_{k-1}+\xi_k+\xi_{k+1}},\quad
k=1,\ldots, N+1.
$$
Therefore, for $k=1,\ldots,N$
\begin{align}
\label{pksym}
 \P\{\zeta_i(t+1)=\zeta_{i}(t)+\delta_{i,k},\ \forall  i \|\zeta(t)=\bx\}
   &\varpropto \beta^{x_{k-1}+x_k+x_{k+1}},\\ \nonumber \\
 \P\{\zeta_{k}(t+1)=\zeta_k(t)-1,\, k=1,\ldots,N \|\zeta(t)=\bx\}
  &\varpropto \beta^{x_{1}+x_{N}},\nonumber
\end{align}
where $x_{N+1}=0$ in (\ref{pksym})  by convention. Let
$$
u_k(t)=\xi_{k-1}(t)+\xi_k(t)+\xi_{k+1}(t),\, k=1,\ldots, N+1,
$$
and recall that $u_k(t)$ is a potential of site $k$ at time $t$.
The process $u(t)=(u_1(t),\dots,u_{N+1}(t))$ is a Markov chain
with the transition probabilities given by
$$
\P(E_k(t)|u(t))=\frac{\beta^{u_k(t)}}{\sum_{i=1}^{N+1} \beta^{u_i(t)
}},\ k=1,2,\dots,N+1
$$
where
 \bn
 E_k(t)=\{&u_{i}(t+1)=u_{i}(t)+1,\ i=k-1,k,k+1,\\
 &\text{ and }  u_{i}(t+1)=u_{i}(t)\ \forall i\notin\{k-1,k,k+1\}\},
\en that is, at time $t+1$ a particle is adsorbed at node $k$.

Fix some small $\eps>0$. First, we will show that there is a
$\delta=\delta(\b,N,\eps)$ such that if for some time $T$ we have
$u_k(T)=\max_i u_i(T)$, then
$$
\P\left( B_{\infty}^{(T,k)} \| \F_T\right)>\delta, \text{ where }
B_{\infty}^{(T,k)}=\bigcap_{t=T}^{\infty} \left[E_{k-1}(t)\cup
E_{k}(t)\cup E_{k+1}(t)\right]
$$
and $\F_T$ is the sigma-algebra generated by the process $u(t)$ up
to time $T$. Secondly, using the conditional Borel-Cantelli lemma,
we will establish that, in fact, with probability $1$ there will
be a $k\in\{1,2,\dots,N+1\}$ for which the event
$B_{\infty}^{(T,k)}$ occurs. Finally, we will show that
$B_{\infty}^{(T,k)}$ implies, with probability one, that either
$\bigcap_{t=T}^{\infty} \left[E_{k-1}(t)\cup E_{k}(t)\right]$ or
$\bigcap_{t=T}^{\infty} \left[E_{k}(t)\cup E_{k+1}(t)\right]$
occurs.

Initially, suppose that $N+1\ge 5$. Without the loss of
generality, assume that $k=3$, i.e., $u_3(T)=\max_i u_i(T)$. Since
the process $u$ is time-homogeneous, we can also set $T=0$. Denote
 \bn
A_{t+1}&=E_2(t+1) \cup E_3(t+1)\cup E_4(t+1),
\\
B_{t}&=\bigcap_{s=1}^{t} A_{s},
\\
 C_{t}&=\left\{\sum_{s=1}^t \left[ 1_{E_2(s)}+1_{E_4(s)}
\right] \le \left(\frac 23+\eps\right) t\right\}
\\
 &=
 \{\text{nodes $2$ and $4$ together adsorb less than $\left( 2/3+\eps\right)t$}
 \\
 &\ {}\ {\ } \ {\ }\text{ particles during the fist $t$ trials}\},
 \en
where $1_E$ is the indicator function of event $E$. In these
notations
\begin{equation}
\label{B_i}
 \P(B_{t+1})\ge \P(A_{t+1}\|B_{t}C_{t})\P(C_{t}\|B_{t})\P(B_{t}).
\end{equation}
Note that $u_3(t)=u_3(0)+t$ on event $B_{t}$. Also
$$
\max\{u_1(t),u_5(t)\}\le u_3(0)+(2/3+\eps) t
$$
on event $B_{t} C_{t} $.  Consequently, we can bound the first
conditional probability in the right side of (\ref{B_i}) as
follows
\begin{align}
\label{A_i} \P(A_{t+1} \| B_{t}C_{t})&=&\nonumber
 \left(1+\frac{[\b^{u_1(t)}+\b^{u_5(t)}]+\sum_{i\notin\{1,\dots,5\}}
 \b^{u_i(t)}} {\b^{u_2(t)}+\b^{u_3(t)}+\b^{u_4(t)}} \right)^{-1}
  \\&\geq& \frac{1}{1+2\beta^{- \left(\frac 13-\eps\right)t} +
(N-4)\beta^{-t}}.
\end{align}
Also, on event $B_{t}$ for $s\le t$ we have $u_2(s)\le u_3(s)$ and
$u_4(s)\le u_3(s)$, hence the probability to get adsorbed at nodes
$2$ or $4$ is smaller than or equal to $2/3$, hence, using large
deviation estimate for the sum of Bernoulli($p$) random variables
from~\cite{Shyr}, Section IV.5, with $p=2/3$, we obtain
\begin{equation}
\label{C_i} \P(C_{t} \| B_{t})\geq 1-2 e^{- 2t\eps^2}.
\end{equation}
Combining bounds (\ref{A_i}) and (\ref{C_i}) we conclude that
$$\P(B_{t+1})\geq (1-\gamma_t)\P(B_{t}),$$
where $\gamma_t$ is summable. Thus event
$B_{\infty}=\cap_{t=1}^{\infty}B_{t}$ occurs with probability at
least
$$
\delta=\P(B_1)\prod_{t=1}^{\infty} (1-\gamma_t)>0.
$$
Let $k^*(T)=\min\{i\in\{1,2,\dots,N+1\}:\ u_i(T)=\max u_j(T)\}$.
Since $\P\left(B_{k^*(T),\infty}^{(T)})\|\F_T\right)\ge \delta$,
by the second conditional Borel-Cantelli lemma, eventually one of
these events will occur.

However, on event $B_{\infty}^{(T,k)}$ it is easy to see that
there will be infinitely many arrivals to the set of nodes
$\{k-1,k+1\}$ (as it is impossible a.s.\ that only node $k$
adsorbs all the particles). Hence, if we consider the times when a
particle is adsorbed at either node $k-1$ or node $k+1$, then
$(u_{k-1},u_{k+1})$ has the same distribution as the balls in a
two-colour urn with exponential reinforcement $w(x)=\b^x$ (see
Section~\ref{stab}); hence by Rubin's theorem we eventually stop
picking one of the two nodes.

Finally, suppose that $\prod_{t=T}^{\infty} (E_{k-1}(t)\cap
E_{k}(t))$ occurred. On this event, for each $t\ge T$, the
probability to get adsorbed at $k$, divided by the probability of
that of $k-1$, is constant and equals $\beta^c$ where
$c={u_{k}(T)-u_{k-1}(T)}=\xi_{k+1}(T)-\xi_{k-2}(T)$, from which
the statement of the Theorem follows by the strong law of large
numbers applied to a sequence of i.i.d.\
Bernoulli($\beta^c/(1+\beta^c)$) random variables and the
finiteness of $u_{k-1}(T)$ and $u_{k+1}(T)$.

\vskip 5mm To finish the proof, we need to consider the special
cases: $N+1=4$. In this case  the proof is virtually identical to
the case $N+1\ge 5$, except that we get slightly different
expression in the RHS of (\ref{A_i}), given by $
{1}/\left[{1+\beta^{- \left(1/3-\eps\right)t} }\right]. $
 \Cox

\begin{rema}  It is relatively easy to prove just transience of the
process of differences.
Namely,  consider  a stochastic process formed by differences of potentials
$$v_k(t)=u_k(t)-u_{N+1}(t),\, k=1,\ldots, N.$$
It is easy to check that process  $(v(t),\, t\in \Z_{+})$ is a Markov chain and that
transience of $(v(t),\, t\in \Z_{+})$  yields transience of
$(\zeta(t),\, t\in \Z_{+})$.
In  turn transience of  process  $(v(t),\, t\in \Z_{+})$ can be established
by applying Theorem \ref{th_fa2} with set
\begin{equation*}
M={\cal C}_a=\left\{
\by\in \Z^N: y_{1}>a,\, y_2-y_k>a,\, k=4,\ldots,N\right\}
\end{equation*}
and Lyapunov function
$$f(\by)=\beta^{-y_{1}}+\sum\limits_{k=4}^{N}\beta^{-y_2+y_k},$$
if $N\geq 4$, and
set $M={\cal D}_a=\left\{
\bx\in \Z^3: y_{2}>a\right\}$ and  Lyapunov  function
$f(\by)=\beta^{-y_{2}},$ if $N=3$,
where in both cases $a>1$ can be  any integer.
\end{rema}

\subsubsection{Proof of Theorem \ref{sym} for $0<\beta<1$}
Assume now  that $0<\beta<1$.
First, assume that $N+1=2M$ is even.
If the process $\zeta$ were recurrent, there would be infinitely
many times $t$ when $\xi_1(t)=\xi_2(t)=\dots=\xi_{N+1}(t)$.
However, we will show that given such a configuration occurs at
time $T$, with a positive probability, independent of $T$, the
following event
 \bn
 A:=\{\lim_{t\to\infty} \xi_i(t)/t=1/M\mbox{ for even $i$ and } \sup_{t\ge T} \xi_i(t)=\xi_i(T)\mbox{ for odd $i$}\}
 \en
occurs. This immediately implies  transience of $\zeta$.

Intuitively, the reason why $A$ occurs with a positive
probability, is the following. Without loss of generality assume
$T=0$ and $\xi_i(T)=0$ for all $i$. Then as long as $\xi_i(t)$
remain $0$ for odd $i$'s for all $t$, the process on even $i$'s is
a \Polya urn scheme, with negative reinforcement, hence we must
have that the relative heights of the ``peaks'' converge to one
(see the proof of Theorem~\ref{nointeraction}). Therefore, for
large times $t\ge t_0$ we have $\xi_i(t)\approx t/M$ for even
$i$'s. On this event, the probability never to add anything into
odd $i$ is asymptotically bounded below by
$$
\prod_{t=t_0}^{\infty}
\left[\frac{M\beta^{t/M}}{M\beta^{t/M}+M\beta^{2t/M}}\right]=\prod_{t=t_0}^{\infty}
\left[1-\frac{1}{1+(1/\beta)^{t/M}}\right]>0
$$
and hence is ``compatible'' with our initial assumption that
$\xi_i$'s remain unchanged for all odd $i$'s.

To make the argument above rigorous, we borrow the idea from the
proof of the main theorem in~\cite{V2001}. Namely, consider the
process $\xi(t)$ at the stopping times $\tau_k$ when the maximum
value of $\xi_{2j}$ reaches $k^2$, that is
$$
\tau_k=\min\{t>0: \max_{j=1,\dots,M} \xi_{2j}(t)=k^2\}, \ \
k=1,2,\dots.
$$
Let event $A_k$ be
$$
A_k=\left\{\min_{j=1,\dots,M} \xi_{2j}(\tau_k) \ge  k^2- k \mbox{
and } \xi_{2j-1}(\tau_k)=0, \ j=1,\dots,M \right\}.
$$
We will show
$$
\P(A_{k+1}\| A_k,A_{k-1},\dots)\ge 1-\gamma_k
$$
with $\sum \gamma_k<\infty$ and hence $\prod_{k=k_0}^{\infty}
(1-\gamma_k)>0$. Since $A\supseteq \bigcap_{k=k_0}^{\infty} A_k$,
this yields
 \bn
 \P(A)&\ge \P\left(\bigcap_{k=k_0}^{\infty}A_{k}\right)
 =\P(A_{k_0})\prod_{k=k_0}^{\infty}\P(A_{k+1}\|A_k,A_{k-1},\dots,A_{k_0})\\
 &\ge  \P(A_{k_0})\prod_{k=k_0}^{\infty} (1-\gamma_k)>0.
 \en
To establish this, first of all, observe that given $A_k$ the
conditional probability that none of odd $\xi_i$ increases between
times $\tau_k$ and $\tau_{k+1}$, is bounded below by
$$
\left(1-\frac{\b^{2k^2-2k}}{\b^{(k+1)^2}}\right)^{\tau_{k+1}-\tau_k}\ge
\left(1-\b^{k^2-4k-1}\right)^{8Mk^2}\ge 1-8Mk^2\b^{k^2-4k-1}
=:1-\gamma'_k
$$
 since
$$
\tau_{k+1}-\tau_{k}\le \tau_{k+1}\le  2M(k+1)^2\le 8Mk^2.
$$
for $k\ge 1$.

From now on let us assume that we are on the event where all
$\xi_{2j-1}$'s remain $0$ during the time interval
$[\tau_k,\tau_{k+1})$; this effectively restricts all adsorptions
to take place at even nodes and thus the transition probabilities
are well-defined. Let
$$
\kappa=\kappa(s):=\min\{t>s: \max_{j=1,\dots,M}
\xi_{2j}(t)=1+\max_{j=1,\dots,M} \xi_{2j}(s)\}
$$
be the time when the maximum increases for the first time after
time $s$.

\begin{lemma}\label{lem_catchup}
There is an $\eps>0$ independent of anything but $M$, such that
$$
\P\left(\mbox{all }\xi_{2j}(\kappa-1) \mbox{ are the same} \|\F_s
\right)\ge \eps,
$$
where $\kappa=\kappa(s)$.
\end{lemma}

Loosely speaking, this states that with probability at least
$\eps$ all $\xi_{2j}$'s will ``catch up'' with the largest of them
before this largest one increases even by $1$, independently of
the past and the actual values of $\xi_{2j}$'s.

\vskip 5mm

 {\sf Proof of Lemma~\ref{lem_catchup}.} Fix a large
$a>0$ such that $\b^a<1/M^2$ and first assume that for some $t\in
[s,\kappa)$
$$
 \max_{j=1,\dots,M} \xi_{2j}(t)-\min_{j=1,\dots,M} \xi_{2j}(t)=l\ge a.
$$
Then, since there are at most $M-1$ maxima, the probability that
each of the minima (and there are at most $M-1$ of them) increases
by $1$ before any of the maxima increases, is at least
$$
 \left[1-(M-1)\b^l\right]^{M-1}>1-M^2\b^l.
$$
However, when all the minima increase, it results in $l\mapsto
l-1$. Consequently, since $1-M^2\b^a>0$, the probability that all
of $\xi_{2j}$ increase until they fall within the $a$-distance of
the maximum before the maximum even changes by one is bounded
below by
$$
\eps'=\prod_{l=a}^{\infty}\left(1-M^2\b^l\right)>0.
$$
But once  $\max_j \xi_{2j}(t)-\min_j \xi_{2j}(t)\le a$, the
probability that all $\xi_{2j}$'s will catch up the maximum value
before $\kappa$ is at least
$$
\eps''=M^{-a(M-1)}>0.
$$
Hence the statement of the Lemma follows with $\eps:=\eps'\times
\eps''$.
 \Cox

\vskip 5mm

 Now pick an $\eps>0$ satisfying the condition of
Lemma~\ref{lem_catchup} and let
$$
\tilde\tau(y)=\min\{t: \ \max_{j=1,\dots,M} \xi_{2j}(t)\ge y\}.
$$
Observe that process $\xi_{2j}(t)$, $j=1,2,\dots,M$ at the times
$\tilde\tau(k^2+k+x)$ for $x=1,2,\dots, k$. By
Lemma~\ref{lem_catchup}, for each of those stopping times
$\tilde\tau_k$ with probability at least $\eps$, {\it
independently of the past}, all $\xi_{2j}$'s become equal before
the time $\tau_{k+1}$ is reached. Hence, the probability that at
least one of these events (all $\xi_{2j}$ are equal) occurs is at
least
$$
1- (1-\eps)^k =1-\gamma''_k.
$$
However, the event $\{$for some $t<\tau_{k+1}$ we have
$\xi_{2j}(t)=k^2+k+x$ for all $j\}$  yields $\min_j
\xi_{2j}(\tau_{k+1})\ge k^2+k=(k+1)^2-(k+1)$, thus implying the
event $A_{k+1}$.

Consequently, we have shown that $ \P(A_{k+1}\| A_k,
A_{k-1},\dots)\ge 1-\gamma_k $ where
$\gamma_k=\gamma'_k+\gamma''_k$ is obviously summable, and this
proves the required transience for the case when $N+1$ is even.

\vskip 5mm

The case of odd $N+1=2M+1$ is analyzed similarly, though obviously
we cannot produce a transient configuration with alternating peaks
and zeros. Instead, redefine the event $A$ as
 \bn
 A:=\left\{\lim_{t\to\infty} \frac{\xi_{1}(t)}t=\lim_{t\to\infty} \frac{\xi_{2}(t)}t=\frac
 1{2M}, \
  \lim_{t\to\infty} \frac{\xi_{i}(t)}t=\frac 1M\mbox{ for even $i\ge 4$}
  \right. \\ \left. \mbox{and }\sup_{t\ge T} \xi_i(t)=\xi_i(T)\mbox{ for odd $i\ge 3$}\right\},
 \en
the positive probability of which will imply transience of process
$\zeta$.

Set
 \bn
  \bar\xi_{2j}(t)&=\left\{\begin{array}{ll}
  \xi_1(t)+\xi_2(t), & \mbox{ if }j=1;\\
  \xi_{2j}(t), & \mbox{ if }j=2,\dots,M,
  \end{array}
  \right.
 \\
  \tau_k&=\min\{t>\tau_{k-1}: \max_{j=1,\dots,M}
  \bar\xi_{2j}(t)=k^2\}, \ \ k=1,2,\dots.
\en
Let now the events $A_k$ be
  \bn
  A_k=\left\{\min_{j=1,\dots,M} \bar\xi_{2j}(\tau_k) \ge  k^2- k, \
  \left|\xi_1(\tau_k)-\xi_2(\tau_k)\right|<k^{1.6},\right.\\
  \left.\mbox{ and }   \xi_{2j+1}(\tau_k)=0, \ j=1,\dots,M \right\}.
  \en
As before, we want to show $\P(A_{k+1}\| A_k,A_{k-1},\dots)\ge
1-\gamma_k$ with $\sum \gamma_k<\infty$.

Firstly, the lower estimate of the conditional probability (given $A_k$) that none of the
$\xi_i$'s for odd $i\ge 3$ increases before $\tau_{k+1}$ is now
replaced by
$$
\left(1-\frac{\b^{\frac
32\left(k^2-k\right)-k^{1.6}}}{\b^{(k+1)^2}}\right)^{\tau_{k+1}-\tau_k}\ge
\left(1-\b^{k^2/3}\right)^{3Mk}\ge 1-3Mk\b^{k^2/3} =:1-\gamma'_k
$$
for $k\ge \max(M,135)$.

Secondly, conditioned on the above event, the probability that
$$\min_{j=1,\dots,M} \bar\xi_{2j}(\tau_{k+1}) \ge  (k+1)^2- (k+1)$$
is bounded below by the same quantity $1-\gamma_k''$ as in the
case when $N$ is even (the proof is a verbatim copy of that case).

Finally, which is new here, we need to ensure that
$$
\left|\xi_1(\tau_{k+1})-\xi_2(\tau_{k+1})\right|<(k+1)^{1.6}=k^{1.6}+1.6
k^{0.6}+o(1).
$$
Now, the conditional probability that $\xi_1$ or $\xi_2$ resp.\
increases, {\it given that either of them increases}, is the same
and equals $1/2$. On the other hand, the number of those times
$t\in[\tau_k,\tau_{k+1}]$ when $\xi_1(t)+\xi_2(t)$ indeed
increases, lies between $k$ and $3k+1$, hence by the law of large
deviations used as in~\cite{V2001}, we have
 \bn
 \P(\left|\left\{\xi_1(\tau_{k+1})-\xi_2(\tau_{k+1})\right\}-
 \left\{\xi_1(\tau_{k})-\xi_2(\tau_{k})\right\}\right|>k^{0.6})\\
 \le \exp(-2k^{0.2}+o(1))=:\gamma_k'''
 \en
Consequently,
$$
\P(A_{k+1}\|A_k, A_{k-1},\dots)\ge 1-\gamma_k
$$
with $\gamma_k=\gamma_k'+\gamma_k''+\gamma_k'''$ and the same
conclusion (i.e.\ transience) holds. \Cox

\section{Criteria}

The following constructive martingale criteria for ergodicity and
transience were used in our paper; their proof can be found
in~\cite{fmm}. Let $(\kappa_k,\, k\geq 0)$ be an irreducible
aperiodic Markov chain on a countable set ${\cal A}$.

\begin{thm}[\cite{fmm}, Theorem 2.2.2]\label{th_fa2}
The Markov chain~$\kappa_k$ is transient if and only if there
exist a set~$M\subset {\cal A}$ and a positive
function~$f(\eta_k)$ such that
$$
  \E \left[ f(\eta_{k+1})-f(\eta_{k} ) \mid \eta_k=v \right]\leq 0
  \ \ \ \ \ \ \  \mbox{for all $v\notin M$;}
$$
$$
  f(v_1)<\inf_{v\in M}f(v),
  \ \ \mbox{
 at least for one
  }v_1\notin M.
$$
\end{thm}
\begin{thm}[\cite{fmm}, Theorem 2.2.3]\label{th_fa1}
 The Markov chain~$\eta_k$ is ergodic if and only if
there exist a finite set~$M$, an $\eps>0$, and a positive
function~$f(\eta_k)$ such that
$$
  \E \left[ f(\eta_{k+1})-f(\eta_{k} ) \mid \eta_k=v \right]\leq -\eps
  \ \ \ \ \ \ \  \mbox{for all $v\notin M$;}
$$
$$
  \E \left[ f(\eta_{k+1}) \mid \eta_k=v \right]<\infty
  \ \ \ \ \ \ \  \mbox{for all $v$.}
$$
\end{thm}

\section*{Acknowledgements}
The authors are grateful to the anonymous referee for meticulous
reading of the paper and pointing out a number of errors /
omissions. We believe that the implementation of his/her
suggestions resulted in a much better paper.

\end{document}